\documentclass[11pt,french]{article}

\usepackage[applemac]{inputenc}
\usepackage{ae, aeguill}
\usepackage[francais]{babel} 
\usepackage{amssymb}
\usepackage{graphicx}
\usepackage{oldgerm}

\def \I{{\mathbb I}}

\def \N{{\mathbb N}}

\def \P{{\mathbb P}}

\def \R{{\mathbb R}}

\def \Z{{\mathbb Z}}

\def\vareps{\varepsilon}

\def \ER{{\cal E}}
\def \FR{{\cal F}}

\def \LR{{\cal L}}

\def \SR{{\cal S}}

\def \v#1{\vspace{#1mm}}
\def \n{{\noindent}}

\def \leq{{\; \leqslant \; }}
\def \geq{{\; \geqslant \; }}

\def \e{{\rm e}}

\def \vareps{\varepsilon}
\def \norm#1#2{\| #1 \|_{#2}}

\def \o#1{\overline{#1}}

\def \fdem{\hfill$\bullet$}

\countdef\annee=110
\countdef\jour=111
\countdef\aux=112
\countdef\bissect=113
\countdef\nbbissect=114

\aux=\year
\divide\aux by 400
\multiply\aux by 400
\annee=\year
\advance \annee by -\aux  % rang de l'année dans le cycle de 400 ans
\aux=\annee
\bissect=\annee
\divide\aux by 4
\nbbissect=\aux
\multiply\aux by 4
\advance \bissect by -\aux % rang de l'année dans le cycle de 4 années
\advance \nbbissect by 1
\ifnum \bissect=0  \advance \nbbissect by -1 \fi   % nombre d'années bissectiles antérieures à l'année en cours

\jour=\annee
\advance\jour by \nbbissect
\ifnum \nbbissect>25 \advance\jour by -1 \fi
\ifnum \nbbissect>50\advance\jour by -1 \fi
\ifnum \nbbissect>75 \advance\jour by -1 \fi

\ifnum \month>1 \advance\jour by 3 \fi
\ifnum \month>3 \advance\jour by 3 \fi
\ifnum \month>4 \advance\jour by 2 \fi
\ifnum \month>5 \advance\jour by 3 \fi
\ifnum \month>6 \advance\jour by 2 \fi
\ifnum \month>7 \advance\jour by 3 \fi
\ifnum \month>8 \advance\jour by 3 \fi
\ifnum \month>9 \advance\jour by 2 \fi
\ifnum \month>10 \advance\jour by 3 \fi
\ifnum \month>11 \advance\jour by 2 \fi

\ifnum \nbbissect=25 \aux=0 \fi
\ifnum \nbbissect=50 \aux=0 \fi
\ifnum \nbbissect=75 \aux=0 \fi

\ifnum \bissect=0 \ifnum \month>2 \advance \jour by \aux \fi\fi % gestion du 29 février éventuel de l'année en cours

\advance\jour by \day

\aux=\jour
\divide\aux by 7
\multiply\aux by 7
\advance \jour by -\aux %rang du jour dans la semaine
\ifnum \jour=0 \jour=7 \fi

 \def\date{ \ifcase\jour\or Samedi \or Dimanche \or Lundi \or Mardi \or Mercredi \or Jeudi \or Vendredi \fi \the\day\ifnum \day=1 $^{\rm er}$ \fi
                         \ \ifcase\month\or Janvier \or Février \or Mars \or Avril \or Mai \or Juin\or Juillet\or Août\or Septembre \or Octobre \or Novembre \or Décembre \fi { \the\year}}

\begin{document}
\centerline{\sc Un théorème limite central local}

\centerline{\sc  en environnement aléatoire stationnaire de conductances sur $\Z$}

\v1
\centerline{\sc A local central limit theorem}

\centerline{\sc  in stationary random environment of conductances on $\Z$}

\v1
\centerline{
Jean-Marc Derrien\footnote[1]{Université de Brest, CNRS - UMR 6205, Laboratoire de Mathématiques de Bretagne Atlantique - 6, avenue Le Gorgeu,
CS 93837,
29238 BREST cedex 3,
France}
}

\v1
\centerline{\it \date}

\v5
\n{\bf Résumé. }
On démontre un théorème limite central local
pour les marches aléatoires aux plus proches voisins
en environnement aléatoire stationnaire de conductances sur $\Z$ en s'affranchissant simultanément des
deux hypothèses classiques d'uniforme ellipticité et d'indépendance sur les conductances.
Outre le théorème limite central, on utilise pour cela des inégalités différentielles discrètes
du type \flqq inégalités de Nash\frqq\ associées
à  la représentation de Hausdorff des suites complètement décroissantes.
La méthode s'adapte aux chaînes de Markov analogues en temps continu.

\v3
\n{\bf Abstract. }We prove a local central limit theorem for nearest neighbours random walks in stationary random environment of conductances on $\Z$ without using any of both classic assumptions of uniform ellipticity and  independence on the conductances. Besides the central limit theorem, we use
discrete differential Nash-type inequalities
associated
with the Hausdorff's representation of the completely decreasing sequences.
The method is also valid for analogous continuous time Markov chains.

\v3
\n{\it 2010 MSC: }60J10; 60K37

\v3
\n{\it Mots-clés: }Marches aléatoires;
environnement aléatoire stationnaire de conductances; théorème limite central local;
inégalités de Nash;
représentation de Hausdorff des suites complètement décroissantes;
théorèmes ergodiques.

\v3
\n{\it Keywords: }Random walks;
stationary random environment of conductances;
local central limit theorem;
Nash's inequalities;
Hausdorff's representation of the completely decreasing sequences;
ergodic theorems.

\v5
\noindent{\bf 1. Le modèle et la présentation des principaux résultats }

\v2
Depuis les années 1980, une part importante de l'étude des milieux aléatoires
a trait aux marches aléatoires aux plus proches voisins en environnement aléatoire stationnaire de conductances
sur $\Z^d$ (dont notamment les marches aléatoires sur l'amas de percolation).
Dans ce cadre particulier,
les différences entre les milieux considérés s'articulent essentiellement
autour des trois modalités suivantes:

{\leftskip=1.5cm
\n - la dimension $d$,

\n - le domaine des valeurs autorisées pour les conductances,

\n - considérer une famille de conductances aléatoires stationnaire ergodique
({\it environnement stationnaire})
vs. considérer le cas particulier d'une famille de conductances indépendantes
et identiquement distribuées ({\it environnement i.i.d.}). 
\par}

\n Par exemple, l'hypothèse d'uniforme ellipticité, c'est-à-dire d'encadrement des conductances par deux constantes strictement positives, permet l'utilisation de nombreuses méthodes classiques dans l'étude des équations aux dérivées partielles.

\v3
Pour faire un historique succinct des travaux qui ont motivé la présente étude, on peut citer, entre autres,
\cite{K} pour une vue d'ensemble des modèles de marches aléatoires en milieux aléatoires sur $\Z^d$
et la systématisation des techniques de martingales et
d'\flqq environnement vu de la particule\frqq\ dans ce contexte, 
\cite{G}, \cite{B} et \cite{B-De} pour l'étude de la récurrence sur l'amas de percolation et en environnement
(de conductances) stationnaire non uniformément elliptique en dimension $d\geq 2$,  
\cite{D-D} pour une étude de la variance asymptotique en environnement stationnaire non uniformément elliptique sur $\Z$,  
\cite{B-D}, \cite{TLC1} et \cite{De} pour un théorème limite central en environnement stationnaire uniformément elliptique
en dimension $d\geq 1$, 
\cite{TLC1}, \cite{TLC2}, \cite{TLC3}, \cite{TLC4}, \cite{TLC5}, \cite{TLC6}, \cite{Biskup}, \cite{TLC7}
et \cite{ADS} pour des théorèmes limites centraux et des principes d'invariance en environnements stationnaires non uniformément elliptiques en dimension $d\geq 2$, 
\cite{Del} pour des inégalités gaussiennes en environnement déterministe uniformément elliptique
en dimension $d\geq 1$, 
\cite{IG1}, \cite{IG2} et \cite{IG3} pour des inégalités gaussiennes 
pour l'amas de percolation en dimension $d\geq 2$, 
\cite{NIG1}, \cite{NIG2} et \cite{NIG3} pour des comportements singuliers du noyau de la chaleur
en environnements i.i.d. non uniformément elliptiques en dimension $d\geq 4$, 
\cite{TLL} pour un théorème limite central local pour l'amas de percolation en dimension $d\geq 2$.

\v3
Dans cet article, on démontre un théorème limite central local
dans le cas stationnaire non uniformément elliptique unidimensionnel.
On constate en particulier que
le rapport entre la moyenne des résistances et celle des conductances  joue un rôle
dans le théorème limite central local
alors que
seul le produit de ces moyennes
apparaît dans l'expression de la variance asymptotique de la marche aléatoire.

\v3
Précisons à présent le modèle étudié.

\v1
On considère en premier lieu un système dynamique probabilisé ergodique
$(\Omega,\FR,\mu,T)$, c'est-à-dire une
mesure de probabilité $\mu$ sur un espace mesurable $(\Omega,\FR)$
et une transformation $T$ de $\Omega$,
inversible, bi-mesurable, préservant la probabilité $\mu$ et pour laquelle
les ensembles invariants sont de mesure égale à 0 ou à 1
(voir \cite{P} par exemple).
Alors,
étant donnée
une application mesurable\\
$c:\Omega \rightarrow ]0,+\infty[$,
on appelle
{\it conductance} de l'arête (non orientée) $\{x,x+1\}$ de $\Z$
dans l'{\it environnement} $\omega$
le réel strictement positif $c (x,x+1) (\omega):=c(T^x \omega)$
(on pose aussi $c (x+1,x) (\omega):=c (x,x+1) (\omega)$).
De la sorte,
la famille de conductances $(c(x,x+1))_{x\in \Z}$
constitue une suite de variables aléatoires stationnaire ergodique.

\v1
Un environnement $\omega$ dans $\Omega$ étant fixé,
on s'intéresse au comportement de la chaîne de Markov
$(S_n)_{n\geq 0}$ sur $\Z$ partant de $x_0\in \Z$
dont les probabilités de transition sont données
par: 
$$
\P_{x_{0}}^\omega [S_{n+1}=x+1\; |\ S_{n}=x]={{c (x,x+1) (\omega)}\over{\o{c} (x) (\omega)}}=:p(x,x+1)(\omega)
$$
et
$$
\P_{x_0}^\omega [S_{n+1}=x-1\; |\ S_{n}=x]={{c (x,x-1)(\omega)}\over{\o{c} (x)(\omega)}}=:p(x,x-1)(\omega)\; ,
$$
où l'on a posé
$\o{c} (x) (\omega):= c(x-1,x) (\omega)+c(x,x+1) (\omega)$, ceci
pour tout $x$ dans $\Z$
($\P^{\omega}_{x_{0}}$ est la probabilité sous laquelle évolue,
dans l'environnement $\omega$,
la chaîne de Markov $(S_n)_{n\geq 0}$ partant de $x_0$).
Notons que la chaîne de Markov $(S_{n})_{n\geq 0}$ est 
réversible sur $\Z$ au sens où:
$$
\forall n\in \N\; ,\ \forall x,y\in \Z\; ,\quad
{{\P_{x}^{\omega} [S_{n}=y]} \over {\o{c} (y)}}
=
{{\P_{y}^{\omega} [S_{n}=x]} \over {\o{c} (x)}}\; .
$$

\v5
L'objet de cet article est de démontrer les trois résultats suivants
qui
prolongent le travail commencé dans \cite{D}. Ils donnent l'ordre de grandeur
de la probabilité pour la marche aléatoire partant de 0 d'être au temps $n$ en un point
$x$ de $\Z$.

\v2
\n{\bf Théorème 1.1 }On a, pour presque tout $\omega$ dans $\Omega$
et pour tout $x$ dans $2\Z$,
$$
\lim_{n\to +\infty}
\sqrt{2n}\, {{\P_{0}^{\omega} [S_{2n}=x]}\over{\o{c} (x) (\omega)}}
=
\left\{
\begin{array}{cl}
{{1}\over{\sqrt{\pi}}}
\,
\sqrt{{{\int {1 \over c}  \, d\mu}\over{\int \o{c}  \, d\mu}}}
&\quad\mbox{si $\o{c}$ et $1/c$ sont intégrables}
\\
 &  \\
0 &\quad\mbox{si $1/c$ est intégrable et si $\int \o{c}\, d\mu=+\infty$}\\
 &  \\
+\infty &\quad\mbox{si $\int {1 \over c}\, d\mu=+\infty$ et si $\o{c}$ est intégrable}
\end{array}
\right.
\ .
$$

\v2
Pour les valeurs de $x$ de l'ordre de $\sqrt{n}$, on montre le

\v1
\n{\bf Théorème 1.2 }
Si $\o{c}$ et $1/c$ sont intégrables alors,
pour presque tout $\omega$ dans $\Omega$,
pour tout réel $x_0$, on a
$$
\lim_{n\to +\infty}
\sqrt{2n}\, {{\P_{0}^{\omega} \left[ S_{2n}=g_n (x_0)\right]}\over{\o{c} \left(g_n (x_0) \right) (\omega)}}
=
{{1}\over{\sqrt{\pi}}}
\,
\sqrt{{{\int {1 \over c}  \, d\mu}\over{\int \o{c}  \, d\mu}}}
\,
{\rm exp} \left( -{{\int \o{c} d\mu\, \int {1\over c} d\mu}\over{4}}x_0^2\right)
\; ,
$$
où $g_n (x_0)$ désigne un plus proche entier pair de $x_0\sqrt{2n}$.

\v2
Le théorème suivant résume les comportements \flqq uniformes\frqq\ obtenus dans le cadre de cette étude.

\v1
\noindent{\bf Théorème 1.3 }
Si $1/c$ est intégrable alors,
pour presque tout $\omega$ dans $\Omega$, pour tout réel $a>0$,
la suite
$$
\left(
\sqrt{2n}\, \max_{x\in B (0,a\sqrt{2n})\cap 2\Z} {{\P_{0}^{\omega} [S_{2n}=x]}\over{\o{c} (x)}}
\right)_{n\geq 0}
$$
est bornée.

Si, de plus, $\int \o{c}  \, d\mu=+\infty$ alors,
pour presque tout $\omega$ dans $\Omega$, pour tout réel $a>0$, on a
$$
\lim_{n\to +\infty}
\sqrt{2n}
\,
\max_{x\in B (0,a\sqrt{2n})\cap 2\Z}
\left(
{{\P_{0}^{\omega} [S_{2n}=x]}\over{\o{c} (x)}}
\right)
=0\; .
$$

\n(On a noté $B (0,a\sqrt{2n})$ l'ensemble des entiers relatifs $z$ tels que $|z|<a\sqrt{2n}$.)

\v3
\n{\it Remarques }

\v1
{\leftskip=1cm
\n - Dans le théorème 1.1, le résultat correspondant au cas où $1/c$ est intégrable mais $\o{c}$ ne l'est pas est un cas particulier du second point du théorème 1.3.

\v1
\n - La preuve du cas non dégénéré du théorème 1.1 et celle du théorème 1.2 étant très similaires,
on détaillera seulement cette dernière (au paragraphe 4).

\v1
\n - A posteriori (car cela n'apparaît pas explicitement dans les preuves), le théorème 1.1 pour $x=0$
doit pouvoir s'interpréter de la façon suivante:
les grandes valeurs prises par $\o{c}$ constituent des \flqq pièges\frqq\ pour la marche aléatoire
$(S_n)_{n\geq 0}$
qui la \flqq retiennent loin de $0$\frqq\ tandis que les grandes valeurs prises par $1/c$ sont des
\flqq barrières\frqq\ qui, au contraire, la \flqq confinent\frqq\ dans des voisinages de $0$.

\v1
\n - Ces théorèmes constituent des \flqq localisations\frqq\ du théorème limite central
pour la marche aléatoire $(S_n)_{n\geq 0}$. Sauf pour le cas où $1/c$ n'est pas intégrable,
on les démontre à partir d'estimations du noyau de la chaleur associé à $(S_n)_{n\geq 0}$
qui résultent elles-mêmes d'inégalités différentielles discrètes,
suivant une démarche initiée par Nash dans \cite{N},
et de
la représentation de Hausdorff des suites complètement décroissantes
(voir la propriété 2.2 ci-dessous pour une définition des suites complètement décroissantes).
L'utilisation de cette démarche de Nash pour obtenir
des estimations sur les probabilités de transition des processus de Markov
réversibles est maintenant très classique (voir \cite{AN1} et \cite{AN2} par exemple).
La représentation de Hausdorff des suites complètement décroissantes
comme suites de moments de lois sur $[0,1]$
(cf aussi le \flqq problème des moments\frqq)
est exposée par exemple dans \cite{F} chapitre VII. 
Le théorème limite central pour la marche aléatoire $(S_n)_{n\geq 0}$ est rappelé ci-dessous.

\par}

\v5
\n{\bf Théorème 1.4 }{\it (Théorème limite central (\cite{K}, \cite{D-D}, \cite{L}))}

\v1
{\leftskip=1cm
\n (1) Si $\o{c}$ et $1/c$ sont intégrables alors,
pour presque tout $\omega$ dans $\Omega$,
pour tous $-\infty \leq a<b\leq +\infty$, on a
$$
\lim_{n\to +\infty}
\P_0^\omega \left[ a<{{S_n}\over{\sqrt{n}}} \leq b\right]
=
\int_{a}^b
k_\sigma (z)
\, dz
\, ,
$$
avec
$$
k_\sigma (z)
:=
{1 \over{\sigma\sqrt{2\pi}}}
\,
\exp \left(-{{z^2}\over{2\sigma^2}}\right)
\quad
{\rm et}
\quad
\sigma^2:={{2}\over{\int \o{c}\, d\mu \, \int {1 \over c}\, d\mu }}\; .
$$

\v1
\n (2) Si $\int \o{c}\, d\mu=+\infty$ ou $\int {1 \over c}\, d\mu=+\infty$ alors,
pour presque tout $\omega$ dans $\Omega$,
pour tous $-\infty \leq a<b\leq +\infty$ avec $a$ et $b$ non nuls, on a
$$
\lim_{n\to +\infty}
\P_0^\omega \left[ a<{{S_n}\over{\sqrt{n}}} \leq b\right]
=
\left\{
\begin{array}{cl}
1 & \quad \mbox{si $0\in ]a,b]$}\\
0 & \quad \mbox{sinon}
\end{array}
\right.
\; .
$$
\par}

\v5
\n{\it Quelques remarques et notations supplémentaires }

\v1
{\leftskip=1cm
\n - Pour tout $x_0$ dans $\Z$ et pour tout réel  $r>0$, on note $B (x_0,r)$
l'ensemble des éléments $x$ de $\Z$
tels que $|x-x_0|<r$.

\v1
\n - On désigne par $\FR (\Z)$ l'ensemble des fonctions définies
sur $\Z$, à valeurs dans $\R$ et à support fini.

\v1
\n - Dans la suite (sauf dans le paragraphe 2 où il n'y a pas d'aléa),
tous les raisonnements sont menés presque sûrement en $\omega$ dans $\Omega$.
On omet systématiquement de mentionner les \flqq$\omega$\frqq.

\v1
\n - Pour toute partie finie $A$ de $\Z$,
on désigne par $\o{c} (A)$ le {\it volume de $A$} défini par:
$$
\o{c} (A)
:=
\sum_{x\in A} \o{c} (x)\; .
$$

\v1
\n - L'opérateur de transition $P$ sur $\Z$ de la marche aléatoire $(S_n)_{n\geq 0}$ est donné par:
$$
Pf(x)
=
f(x-1)\, p(x,x-1)+ f(x+1)\, p (x,x+1)\; ,\quad x\in \Z\; ,
\ f\in \FR (\Z)\; .
$$
Il
préserve l'ensemble $\FR (\Z)$.

\v1
\n - On munit $\FR (\Z)$ d'un produit scalaire et de la norme associée
en posant pour $f$ et $g$ dans $\FR (\Z)$:
$$
(f,g):=\sum_{x\in \Z} f(x)\, g(x)\, \o{c} (x)
\quad
{\rm et}
\quad
\norm{f}{}^2
:=(f,f)\; .
$$

La réversibilité de la marche $(S_n)_{n\geq 0}$
s'exprime alors de manière fonctionnelle comme suit:
$$
\forall f,g \in \FR (\Z)\; ,\quad
(f,Pg)=(Pf,g)\; .
$$
En particulier, $(\o{c} (x))_{x\in \Z}$ est une mesure sur $\Z$
invariante sous l'action de $P$.

\par}

\v5
\n{\it Plan de l'article --}
Dans le paragraphe suivant, on étudie l'équation de la chaleur discrète
associée à une chaîne de Markov réversible
et on remarque en particulier la \flqq complète décroissance\frqq\ au
cours du temps
de l'\flqq énergie\frqq\ de la solution d'une telle équation.
La représentation
de Hausdorff des suites complètement décroissantes nous permet alors de majorer
la différence
entre deux valeurs successives de cette \flqq énergie\frqq. Cette majoration
est appliquée à la marche aléatoire $(S_n)_{n\geq 0}$ aux paragraphes 3 et 4.
Deux inégalités gaussiennes ainsi que
les cas dégénérés du théorème 1.1
sont obtenus dans le paragraphe 3.
Les théorèmes 1.2 et 1.3 sont démontrés dans le paragraphe 4 à la suite d'une
étude de la régularité
du noyau de la chaleur associé à $(S_n)_{n\geq 0}$.
Dans le dernier paragraphe, on s'intéresse au cas d'une chaîne de Markov
analogue en temps continu.

\v5
\noindent{\bf 2. Sur l'équation de la chaleur associée à une chaîne de Markov réversible }

\v2
Dans ce paragraphe, on considère une probabilité de transition $Q=(q_{ij})_{i,j\in \SR}$ sur un espace
d'états dénombrable $\SR$ et une mesure positive $\lambda=(\lambda_{i})_{i\in \SR}$ sur $\SR$ que l'on suppose
réversible pour $Q$. On a donc:
$$
\forall i,j\in \SR\; ,
\quad
\lambda_i\, q_{ij}=\lambda_j\, q_{ji}\; .
$$

On munit l'ensemble $\FR (\SR)$ des fonctions définies sur $\SR$, à valeurs dans $\R$ et à support fini,
d'un produit scalaire et d'une norme en posant:
$$
(f,g):=\sum_{i\in \SR} f(i)\, g(i)\,\lambda_i 
\quad
{\rm et}
\quad
\norm{f}{}^2
:=(f,f)\; .
$$

Pour simplifier la rédaction et compte tenu des applications
en vue, {\bf on suppose pour la suite que, pour tout élément $f$ de $\FR (\SR)$,
$Qf$ est encore élément de $\FR (\SR)$.}

Un calcul élémentaire utilisant la réversibilité de $\lambda$ pour $Q$ montre que:
$$
\forall f,g \in \FR (\SR)\; ,\quad
(f,Qg)=(Qf,g)\; .
$$

\v5
La {\it forme de Dirichlet associée
à $Q^2$ et $\lambda$} est définie par:
$$
\ER(f,g)
:=
{1\over 2}
\sum_{i,j\in \SR}
\left( f(i) - f(j) \right)\, 
\left( g(i) - g(j) \right)\,
q_{ij}^{(2)}\, \lambda_i \; ,
\quad
f,g\in \FR (\SR)\; ,
$$
où l'on a noté
$
q_{ij}^{(2)}
=
\sum_{k\in \SR} q_{ik}\, q_{kj}
$
le coefficient d'indice $(i,j)$ de $Q^2$.

\v3
Le lemme suivant est l'analogue de la première formule de Green pour 
le laplacien discret $I-Q^2$.

\v1
\n{\bf Lemme 2.1: }Pour tout élément $f$ de $\FR (\SR)$, on a:
$$
\ER (f,f)
=
(f,(I-Q^2)f)
=
\norm{f}{}^2
-
\norm{Qf}{}^2
\; .
$$

En particulier, pour tout entier $k\geq 0$, on a:
$$
(f,(I-Q^2)^k f)\geq 0\; .
$$

\vspace{2mm}
\noindent{\it Preuve }L'argument est calculatoire. On a
\begin{eqnarray*}
\ER (f,f)
&=&
{1\over 2}
\sum_{i,j\in \SR}
\left( f(i) - f(j) \right)^2\, 
q_{ij}^{(2)}\, \lambda_i
\\
&=&
{1\over 2}
\sum_{i\in \SR}
f (i)^{2}
\,
\left(
\sum_{j\in \SR}
q_{ij}^{(2)}
\right)
\,
\lambda_i
+
{1\over 2}
\sum_{j\in \SR}
f (j)^{2}
\,
\left(
\sum_{i\in \SR}
q_{ij}^{(2)}\, \lambda_i
\right)
-
\sum_{i,j\in \SR}
f(i) f(j)
\, 
q_{ij}^{(2)}\, \lambda_i
\\
&=&
\sum_{i\in \SR}
f (i)^{2}
\,
\lambda_i
-
\sum_{i,j\in \SR}
f(i) f(j)
\, 
q_{ij}^{(2)}\, \lambda_i
\end{eqnarray*}
en utilisant la réversibilité de $\lambda$ pour $Q^2$.
En utilisant  cette fois la réversibilité de $\lambda$ pour $Q$,
on obtient donc:
$$
\ER (f,f)
=
\sum_{i\in \SR}
f (i)
\left(
f(i)
-
\sum_{j\in \SR}
f(j)
\, 
q_{ij}^{(2)}
\right)
\, \lambda_i
=
(f,(I-Q^2)f)
=
\norm{f}{}^2
-
\norm{Qf}{}^2
\; .
$$

\v1
Ainsi, toujours par réversibilité de $\lambda$, il vient, pour $k=2l$ pair,
$$
(f,(I-Q^2)^k f)=((I-Q^2)^l f,(I-Q^2)^l f)\geq 0
$$
et, pour  $k=2l+1$ impair,
$$
(f,(I-Q^2)^k f)=((I-Q^2)^l f,(I-Q^2)\, ((I-Q^2)^l f))\geq 0\; .
$$

\fdem

\v5
Etant donné $f$ dans $\FR (\SR)$, on considère à présent la suite $(f_n)_{n\geq 0}$
d'éléments de $\FR (\SR)$ solution de
l'{\it équation de la chaleur discrète} (associée à $Q$ et $f$):
$$
\left\{
\begin{array}{l}
f_0=f\; ,\\
f_n - f_{n+1}=(I-Q)f_n\; ,\quad n\geq 0\; .
\end{array}
\right.
$$
On a donc, pour tout $n\geq 0$,
$f_n=Q^n f$.

Intéressons-nous plus particulièrement
à la suite des\ \flqq énergies\frqq\ $\left(\norm{Q^n f}{}^2\right)_{n\geq 0}$
d'une telle solution
et, pour cela, introduisons
la famille des suites de ses \flqq différences successives\frqq\ en notant,
pour tout $k\geq 0$,
$$
\Delta_n^{(k+1)} (f):=\Delta_n^{(k)} (f)- \Delta_{n+1}^{(k)} (f)\;  ,\quad n\geq 0\; ,
$$
où l'on a posé $\Delta_n^{(0)} (f):=\norm{Q^n f}{}^2$ pour tout $n\geq 0$.

\v3
La propriété suivante est élémentaire mais elle joue un rôle déterminant dans la suite.

\v1
\n{\bf Propriété 2.2 }Pour tout élément $f$ de $\FR (\SR)$,
la suite $(\norm{Q^n f}{}^2)_{n\geq 0}$ est {\it complètement décroissante} au sens où:
$$
\forall k,n\in \N\, ,
\quad
\Delta_n^{(k)} (f)\geq 0
\; .
$$

\v3
\n{\it Preuve }C'est une conséquence du lemme 2.1
une fois vérifié par récurrence sur $k$ et en utilisant la
réversibilité de $\lambda$ que
$$
\forall k,n\in \N\, ,
\quad
\Delta_n^{(k)} (f)=\left(Q^n f,(I-Q^2)^k Q^n f\right)\; .
$$

\fdem

\v3
\n{\bf Corollaire 2.3 }
Pour tout élément $f$ de $\FR (\SR)$,
il existe une mesure de probabilité borélienne $\nu_f$ sur $[0,1]$ telle que:
$$
\forall n\geq 0\; ,
\quad
\norm{Q^n f}{}^2=\norm{f}{}^2\, \int_{[0,1]} x^n \, d\nu_f (x)\; .
$$

\n{\it Preuve }C'est une application directe du théorème de Hausdorff
de représentation des suites complètement décroissantes
(voir par exemple \cite{F}\ chapitre VII).

\fdem

\newpage
\n{\it Remarques }

\v1
{\leftskip=1cm
\n - Le  théorème de Hausdorff
de représentation des suites complètement décroissantes (\cite{H})
est à rapprocher du
théorème d'Herglotz de représentation des suites de type positif (\cite{He});
la première représentation est donnée par la suite des moments d'une loi sur [0,1],
la seconde par la transformée de Fourier d'une telle loi.

\v1
\n - La théorie spectrale des opérateurs symétriques
dans les espaces de Hilbert
permet d'obtenir cette même représentation
de la suite $\left( \norm{Q^n f}{}^2 \right)$
mais de manière moins élémentaire
(voir \cite{R} par exemple).

\par}

\v5
\n{\bf Corollaire 2.4 }
Pour tout élément $f$ de $\FR (\SR)$,
pour tout entier $n\geq 0$,
on a:
$$
\norm{Q^{2n} f}{}^2 - \norm{Q^{2n+1} f}{}^2
\leq
{{1}\over {n}}
\,
\norm{Q^{n} f}{}^2\; .
$$

\v2
\n{\it Preuve }
En utilisant la représentation de Hausdorff précédente, il vient:
\begin{eqnarray*}
\norm{Q^{2n} f}{}^2 - \norm{Q^{2n+1} f}{}^2
&=&\norm{f}{}^2\, \int_{[0,1]} x^{2n} (1-x)\, d\nu_f (x)\\
&\leq&\max_{x\in [0,1]}
\left( x^n (1-x) \right)\, \left(\norm{f}{}^2\, \int_{[0,1]} x^{n} \, d\nu_f (x)\right) \\
&=&
{{n^n}\over {(n+1)^{n+1}}}
\,
\norm{Q^n f}{}^2
\leq
{{1}\over {n}}
\,
\norm{Q^n f}{}^2
\; .
\end{eqnarray*}

\fdem

\v5
On suppose à présent que, pour tout $i$ dans $S$, $\lambda_i$ est strictement positif.
De plus,
on particularise un élément $o$ de $\SR$.

On définit
le {\it noyau de la chaleur} $(h_n)_{n\geq 0}$ (\flqq basé en $o$ \frqq\ et)
associé à la probabilité de transition $Q$
comme étant la solution de l'équation de la chaleur discrète associée à $Q$
et à la fonction $h_0$ donnée par:
$$
h_{0} (i)={1 \over{\lambda_o}}\, \I_{\{o\}} (i)
=
\left\{
\begin{array}{cl}
{1 \over{\lambda_o}} & \quad \mbox{si $i=o$}\\
0 & \quad \mbox{sinon}
\end{array}
\right.
\; , \quad i\in \SR\; .
$$

Il vient, pour tout entier $n\geq 0$ et pour tout $i$ dans $\SR$,
$$
h_{n} (i)
=
Q^n h_{0} (i)
=
{{\P_{i} [X_{n}=o]} \over {\lambda_o}}
=
{{\P_{o} [X_{n}=i]} \over {\lambda_i}}\; ,
$$
où $(X_n)_{n\in \N}$ désigne la chaîne de Markov
sur $\SR$ de probabilité de transition $Q$ (partant de $i$ sous la probabilité $\P_{i}$).

\v3
Le lemme suivant est très classique. Il permet notamment de reformuler le théorème 1.1
(pour $x=0$)
en terme d'\flqq énergie\frqq\ du noyau de la chaleur.

\v1
\n{\bf Lemme 2.5 }Pour tout entier $n\geq 0$, on a:
$$
\norm{h_n}{}^2
=
{{\P_o [X_{2n}=o]} \over{\lambda_o}}
=
h_{2n} (o)\; .
$$

\v1
\n{\it Preuve }On a, grâce à la réversibilité de $\lambda$,
$$
\norm{h_n}{}^2
=
(Q^n h_0, Q^n h_0)
=
(Q^{2n} h_0, h_0)
=
Q^{2n} h_0 (o)=h_{2n} (o)\; .
$$
\fdem

\v5
\noindent{\bf 3. Des inégalités gaussiennes et les cas dégénérés}

\v2
On revient à présent à la marche aléatoire $(S_n)_{n\geq 0}$ du paragraphe 1.

Tout au long de ce paragraphe 3, 
on omet de préciser que les résultats obtenus ont lieu pour presque tout environnement
de conductances.

Les résultats du paragraphe précédent s'appliquent à l'opérateur $P$
et à la mesure réversible $(\o{c} (x))_{x\in \Z}$.
Dans ce contexte, la forme de Dirichlet associée à $P^2$
est donnée, pour tous $f$ et $g$ dans $\FR (\Z)$, par:
$$
\ER(f,g)
=
\sum_{x\in \Z}
\left( f(x-1) - f(x+1) \right)\, 
\left( g(x-1) - g(x+1) \right)\, 
{{c(x-1,x)\, c(x,x+1)}\over{\o{c} (x)}}\; .
$$
(La quantité
$$
{{c(x-1,x)\, c(x,x+1)}\over{\o{c} (x)}}
=
{1 \over
{
{{1}\over{c(x-1,x)}}+{{1}\over{c(x,x+1)}}
}}
$$
est égale, comme il se doit, à la conductance
d'un circuit électrique constitué de
deux conductances $c(x-1,x)$ et $c(x,x+1)$ disposées en série.)

\v1
On note encore $(h_n)_{n\geq 0}$ le noyau de la chaleur associé à $(S_n)_{n\geq 0}$
partant de 0:
$$
\forall n\geq 0\; ,\quad
h_n=P^n \left( {{1}\over{\o{c} (0)}} \, \I_{\{0\}}\right)\; .
$$

\v3
Commençons par établir la convergence vers 0 de la suite
$\left(\norm{h_n}{}^2\right)_{n\geq 0}$.

\v1
\n{\bf Proposition 3.1 }La marche $(S_n)_{n\geq 0}$ est
une chaîne de Markov récurrente nulle.

On a donc en particulier la convergence:
$$
\lim_{n \to +\infty} \norm{h_n}{}^2=\lim_{n\to +\infty} {{\P_0 [S_{2n}=0]}\over{\o{c} (0)}}=0\; .
$$

\v2
\n{\it Preuve }Pour tout $K$ dans $\N^*$,
la probabilité pour $(S_n)_{n\geq 0}$ d'atteindre $\{-K,K\}$
avant de revenir en 0 est donnée
par
$$
{{1}\over{\o{c} (0)}}
\,
\left(
{{1}\over{\sum_{x=-K}^{-1} {{1}\over{c(x,x+1)}} }}
+
{{1}\over{\sum_{x=0}^{K-1} {{1}\over{c(x,x+1)}} }}
\right)
\; ,
$$
la {\it conductance effective} entre 0 et $\{-K,K\}$ divisée par $\o{c} (0)$
(\cite{D-S}). Or,
le théorème de récurrence de Poincaré (voir \cite{P} par exemple)
garantit les deux divergences:
$$
\sum_{x=-\infty}^{-1} {{1}\over{c(x,x+1)}}=+\infty
\quad{\rm et}\quad
\sum_{x=0}^{+\infty} {{1}\over{c(x,x+1)}}=+\infty
\; .
$$
On en déduit la récurrence de $(S_n)_{n\geq 0}$.
La récurrence nulle de $(S_n)_{n\geq 0}$ résulte alors
de ce que $(\o{c} (x))_{x\in \Z}$ est une mesure sur $\Z$, $P$-invariante et
de masse totale infinie (toujours d'après le théorème de récurrence
de Poincaré).

La convergence de la suite
$\left( \P_0 [S_{2n}=0] \right)_{n\geq 0}$ vers 0 est une conséquence de cette récurrence
nulle comme il est exposé par exemple dans \cite{G-S} p. 214.

\fdem

\v3
Le théorème suivant donne un minorant (asymptotique)
de la suite $\left(\sqrt{n} \, \norm{h_n}{}^2\right)_{n\geq 0}$.

\v1
\noindent{\bf Théorème 3.2 }
Si $\o{c}$ est intégrable alors
$$
\liminf_{n\to +\infty}
\sqrt{n}
\,
\norm{h_n}{}^2\geq {1\over{13}}
\,
\sqrt{{{\int {1 \over c}  \, d\mu}\over{\int \o{c}  \, d\mu}}}\; .
$$

En particulier,
si $\o{c}$ est intégrable
et
si $\int {1 \over c} \, d\mu=+\infty$ alors
$$
\lim_{n\to +\infty}
\sqrt{n}
\,
\norm{h_n}{}^2
=+\infty\; .
$$

\v3
\n{\it Preuve }
On va exploiter le fait que, typiquement, jusqu'à l'instant $n$,
la marche aléatoire $(S_n)_{n\geq 0}$ \flqq évolue\frqq\ dans un voisinage
de 0 dont la longueur est de l'ordre de $\sqrt{n}$.

\v1 
Pour tout $\delta>0$ et pour tout entier $n\geq 1$, on a:
$$
\sqrt{n}\, \norm{h_n}{}^2
\geq
\sqrt{n}\, \sum_{x\in B(0,\delta \sqrt{n})}  h_{n} (x)^2 \; \o{c} (x)
=
\sqrt{n}\, \sum_{x\in B(0,\delta \sqrt{n})}  {{\P_{0} [S_{n}=x]^2}\over{ \o{c} (x)}}
\; .
$$
Or,
en utilisant l'inégalité de Cauchy-Schwarz, il vient:
$$
\left(\sum_{x\in B(0,\delta \sqrt{n})} \P_{0} [S_{n}=x]\right)^{2}
\leq 
\o{c} \left( B(0,\delta \sqrt{n}) \right)
\,
\sum_{x\in B(0,\delta \sqrt{n})}  {{\P_{0} [S_{n}=x]^2}\over{ \o{c} (x)}}
\; .
$$
On en déduit que
\begin{eqnarray*}
\sqrt{n}\, \norm{h_n}{}^2
&\geq&
{{\sqrt{n}}\over{\o{c} \left( B(0,\delta \sqrt{n}) \right)}} 
\,
\left(\sum_{x\in B(0,\delta \sqrt{n})} \P_{0} [S_{n}=x]\right)^{2}
\\
&=&
{1\over{2 \delta}}
\,
{{2 \delta \sqrt{n}}\over{\o{c} \left( B(0,\delta\sqrt{n}) \right)}} 
\,
\P_{0} \left[ \left| {{S_{n}}\over{\sqrt{n}}} \right|<\delta \right]^{2}\; .
\end{eqnarray*}

\v3
On applique à présent le théorème ergodique ponctuel
de Birkhoff à $\o{c}$ (voir \cite{P} par exemple)
et le théorème limite central (théorème 1.4 ci-dessus) en distinguant selon que $1/c$
est intégrable ou non.

\v1
\n - Si $\int {1 \over c}  \, d\mu=+\infty$ alors, pour tout $\delta>0$,
$$
\liminf_{n\to +\infty} \sqrt{n}\, \norm{h_n}{}^2
\geq
{1\over{2 \delta \int \o{c} \, d\mu}}
\; ;
$$
ce qui donne, en passant à la limite $\delta\to 0^+$, la convergence:
$$
\lim_{n\to +\infty} \sqrt{n}\, \norm{h_n}{}^2
=+\infty\; .
$$

\v1
\n - Si $1/c$ est intégrable alors
il vient,
pour tout $\delta>0$
et
avec
les notations du théorème 1.4,
$$
\liminf_{n\to +\infty} \sqrt{n}\; \norm{h_n}{}^2
\geq
{1\over{2\delta}}
\,
{{1}\over{\int \o{c}\, d\mu}} 
\,
\left(
\int_{-\delta}^\delta
k_\sigma (z)
\, dz
\right)^{2}
\; .
$$
Soit encore:
$$
\liminf_{n\to +\infty} \sqrt{n}\; \norm{h_n}{}^2
\geq
{1\over{4\pi}}
\,
{{1}\over{\sigma \int \o{c}  \, d\mu}} 
\,
{{\sigma}\over{\delta}} 
\left(
\int_{-\delta/\sigma}^{\delta/\sigma}
\exp \left(-{{z^2}\over{2}}\right)
\; dz
\right)^{2}
\; .
$$
En prenant par exemple $\delta=\sigma$, on obtient
par une minoration grossière:
$$
\liminf_{n\to +\infty} \sqrt{n}\; \norm{h_n}{}^2
\geq
{{1}\over{\pi \e\sqrt{2}}}
\,
\sqrt{{{\int {1 \over c}  \, d\mu}\over{\int \o{c}  \, d\mu}}}
\geq
{{1}\over{13}}
\,
\sqrt{{{\int {1 \over c}  \, d\mu}\over{\int \o{c}  \, d\mu}}}
\; .
$$

\fdem
 
\v3
\n{\it Remarque }Le théorème 3.2 assure ainsi que, lorsque $\o{c}$ est intégrable sans que $1/c$ ne le soit, on a, pour presque tout environnement et pour tout $x$ dans $2\Z$,
\begin{eqnarray*}
\lim_{n \to +\infty} \sqrt{2n}\, {{\P_{x} [S_{2n}=x]}\over{\o{c} (x)}}=+\infty
\end{eqnarray*}
(on a utilisé également le lemme 2.5 et la stationnarité de l'environnement).
On en déduit la partie du théorème 1.1 correspondante à ces hypothèses sur les conductances
en utilisant l'irréductibilité de $(S_n)_{n\geq 0}$ et la propriété de Markov.

\v5
On s'intéresse à présent à la majoration
de la suite
$\left(\sqrt{n} \, \norm{h_n}{}^2\right)_{n\geq 0}$.

\v1
\noindent{\bf Théorème 3.3 }
Si $1/c$ est intégrable alors
$$
\limsup_{n\to +\infty}
\sqrt{n}
\,
\norm{h_n}{}^2
\leq
4
\,
\sqrt{{{\int {1 \over c}  \, d\mu}\over{\int \o{c}  \, d\mu}}}
\; .
$$

En particulier,
si $1/c$ est intégrable
et
si $\int \o{c}  \, d\mu=+\infty$ alors
$$
\lim_{n\to +\infty}
\sqrt{n}
\,
\norm{h_n}{}^2
=
0\; .
$$

\v3
\n{\it Preuve }
Pour tout $n\geq 1$, pour tout $K\geq 1$, on considère
un élément
$x_0$ dans $B(0,K)\cap 2\Z$ tel que
$$
h_{2n} (x_{0}) = \min \{h_{2n} (x)\; |\ x\in B(0,K)\cap 2\Z \}\; .
$$

Il vient:
$$
h_{2n} (x_{0})
\leq
{1 \over{\sum_{x\in B(0,K)\cap 2\Z} \o{c}(x)}}
\sum_{x\in B(0,K)\cap 2\Z} h_{2n} (x) \o{c} (x)
\leq
{1 \over{\sum_{x\in B(0,K)\cap 2\Z} \o{c}(x)}}
\; .
$$

En supposant par exemple $x_0=2l_0+2$ et $l_0\geq 0$, on obtient donc:
\begin{eqnarray*}
& &
\norm{h_n}{}^2 - {1 \over{\sum_{x\in B(0,K)\cap 2\Z} \o{c}(x)}}
\\
&\leq&
h_{2n} (0) - h_{2n} (x_0)
\\
&\leq&
\sum_{l=0}^{l_{0}}
\left| h_{2n} (2l) - h_{2n} (2l+2) \right|
\\
&\leq&
\sqrt{\sum_{l=0}^{l_{0}}(h_{2n} (2l) - h_{2n} (2l+2))^{2} \;
{{c(2l,2l+1)\; c(2l+1,2l+2)}\over{\o{c}(2l+1)}}}
\\
& &\times\,
\sqrt{
\sum_{l=0}^{l_{0}}
{1\over{c(2l,2l+1)}}
+
{1\over{c(2l+1,2l+2)}}
}
\; ,
\end{eqnarray*}
en utilisant l'inégalité triangulaire puis
l'inégalité de Cauchy-Schwarz.
Ainsi, compte tenu de l'expression de la forme de Dirichlet
associée à $P^2$ et du lemme 2.1, il vient:
\begin{eqnarray*}
& &
\norm{h_n}{}^2 - {1 \over{\sum_{x\in B(0,K)\cap 2\Z} \o{c}(x)}}
\\
&\leq&
\sqrt{\ER \left(h_{2n},h_{2n} \right)}
\;
\left(\sum_{x=-K}^{K-1} {1\over {c(x,x+1)}}\right)^{1/2}
\\
&=&
\sqrt{\norm{h_{2n}}{}^2 - \norm{h_{2n+1}}{}^2}
\;
\left(\sum_{x=-K}^{K-1} {1\over {c(x,x+1)}}\right)^{1/2}
\; .
\end{eqnarray*}
Le corollaire 2.4 assure alors que
\begin{eqnarray}\label{equation1}
\norm{h_n}{}^2 - {1 \over{\sum_{x\in B(0,K)\cap 2\Z} \o{c}(x)}}
\leq
{{1}\over{\sqrt{n}}}
\,
\norm{h_n}{}
\;
\left(\sum_{x=-K}^{K-1} {1\over {c(x,x+1)}}\right)^{1/2}
\end{eqnarray}

\v2
Pour tout entier $n$ assez grand, on définit à présent $K (n)\geq 1$
par la relation
\begin{eqnarray}\label{equation2}
\sum_{x\in B(0,K(n)-1)\cap 2\Z} \o{c}(x)
\leq
{{3}\over{\norm{h_n}{}^2}}
<
\sum_{x\in B(0,K(n))\cap 2\Z} \o{c}(x)
\end{eqnarray}
Remarquons que la proposition 3.1 garantit que
$\lim_{n\to +\infty} K(n)=+\infty$.

\v2
Il vient d'après (\ref{equation1}) et la seconde inégalité de (\ref{equation2}):
\begin{eqnarray*}
{{1}\over{\sqrt{n}}}
\,
\norm{h_n}{}
&\geq&
{2\over 3}\, \norm{h_n}{}^2\,
\left(\sum_{x=-K(n)}^{K(n)-1} {1\over {c(x,x+1)}}\right)^{-1/2}
\\
&=&
{{\sqrt{2}}\over{3}}\,
\norm{h_n}{}^2\,
\left({{1}\over{2K(n)}} \sum_{x=-K(n)}^{K(n)-1} {1\over {c(x,x+1)}}\right)^{-1/2}
\\
& &\times
\left({{\sum_{x\in B(0,K(n)-1)\cap 2\Z} \o{c} (x)}\over{K(n)}}\right)^{1/2}
\,
\left(\sum_{x\in B(0,K(n)-1)\cap 2\Z} \o{c} (x)\right)^{-1/2}
\; ;
\end{eqnarray*}
ce qui donne en utilisant cette fois la première inégalité de (\ref{equation2}):
\begin{eqnarray*}
& &
{{1}\over{\sqrt{n}}}
\,
\norm{h_n}{}
\\
&\geq&
{{\sqrt{2}}\over{3}}\,
\norm{h_n}{}^2\,
\left({{1}\over{2K(n)}} \sum_{x=-K(n)}^{K(n)-1} {1\over {c(x,x+1)}}\right)^{-1/2}
\left({{\sum_{x\in B(0,K(n)-1)\cap 2\Z} \o{c} (x)}\over{K(n)}}\right)^{1/2}
\left({{\norm{h_n}{}^2}\over{3}}\right)^{1/2}
\; .
\end{eqnarray*}

On achève la preuve du théorème 3.3
en divisant par $\norm{h_n}{}/\sqrt{n}$
et
en appliquant à $1/c$ et à $\o{c}$
le théorème ergodique ponctuel de Birkhoff.

\fdem

\newpage
\noindent{\bf 4. Un théorème limite central local non dégénéré }

\v3
On commence par établir une propriété de régularité du noyau de la chaleur.

\v1
\noindent{\bf Théorème 4.1 }
Si $1/c$ est intégrable alors, pour presque tout environnement,
il existe un réel $C>0$ tel que,
pour tout réel $x_0$,
pour tout réel $a>0$,
il existe un entier $n_0$ tel que,
pour tout entier $n\geq n_0$,
on a:
$$
\max_{x\in B \left(g_n (x_0),a \sqrt{2n}\right)\cap 2\Z}
\left| h_{2n} (x) - h_{2n} \left( g_n (x_0) \right) \right|
\leq
C\, {{\sqrt{a}}\over{\sqrt{n}}}\; ,
$$
où $g_n (x_0)$ désigne comme précédemment
un plus proche entier pair de $x_0\sqrt{2n}$.

\v2
\n{\it Preuve }En utilisant le théorème 3.3 ainsi que sa preuve,
on obtient que, pour presque tout environnement, il existe une constante $C'$
telle que, 
pour tout $x_0$, pour tout $a>0$,
pour tout entier $n\geq 1$,
pour tout $x$ dans $B (g_n (x_0),a \sqrt{2n})\cap 2\Z$, on a:
\begin{eqnarray}
\left| h_{2n} (x) - h_{2n} (g_n (x_0)) \right|
&\leq&
{{1}\over{\sqrt{n}}}
\,
\norm{h_{n}}{}
\;
\left(\sum_{z\in B (g_n (x_0),a \sqrt{2n})} {1\over {c(z,z+1)}}\right)^{1/2}
\label{major}
\\
&\leq&{{C'}\over{\sqrt{n}}}\, \left({{1}\over{\sqrt{n}}}\, \sum_{z\in B (g_n (x_0),a \sqrt{2n})} {1\over {c(z,z+1)}}\right)^{1/2}
\label{ineg1}
\end{eqnarray}

Par ailleurs, d'après le théorème ergodique ponctuel de Birkhoff,
pour presque tout $\omega$ dans $\Omega$,
pour tout $x_0$ et tout $a$ avec $0<a<x_0$,
il existe un entier $N_0=N_0 (\omega,x_0,a)\geq 0$ tel que, pour tout $N\geq N_0$, on a:
$$
N \left( \int {{1}\over{c}}\, d\mu - {{a}\over{x_0}} \right)
\leq \sum_{z=0}^{N-1} {{1}\over{c(z,z+1)(\omega)}}
\leq N \left( \int {{1}\over{c}}\, d\mu + {{a}\over{x_0}} \right)
\; .
$$
Ainsi, en considérant un entier $n_0 = n_0 (\omega,x_0,a)$ tel que, pour tout $n\geq n_0$,
on ait l'inégalité $g_n (x_0) - a \sqrt{2n} \geq N_0 \;$,
il vient:
\begin{eqnarray}
& &
\sum_{z\in B (g_n (x_0),a \sqrt{2n})} {1\over {c(z,z+1)(\omega)}}
\nonumber
\\
&\leq&
\sum_{0 \leq z \leq g_n (x_0) + a \sqrt{2n}}  {1\over {c(z,z+1)(\omega)}}
-
\sum_{0 \leq z \leq g_n (x_0) - a \sqrt{2n}}  {1\over {c(z,z+1)(\omega)}}
\nonumber
\\
&\leq&
(g_n (x_0) + a \sqrt{2n}) \left( \int {{1}\over{c}}\, d\mu + {{a}\over{x_0}} \right)
-
(g_n (x_0) - a \sqrt{2n} - 1) \left( \int {{1}\over{c}}\, d\mu - {{a}\over{x_0}} \right)
\nonumber
\\
&\leq&
a\,\sqrt{2n}\left(5+3 \int {{1}\over{c}} \, d\mu \right)
\label{ineg2}
\end{eqnarray}
(on a supposé pour simplifier $n$ assez grand pour avoir l'inégalité $a\sqrt{2n} \geq 1$).

\v1
Les inégalités (\ref{ineg1}) et (\ref{ineg2}) permettent de conclure dans le cas $0<a<x_0$.
Les autres cas se traitent de manière analogue.

\fdem

\newpage
\n{\it Remarques }

\v1
{\leftskip=1cm
\n - Dans l'approche classique de Moser (\cite{M}),
un tel résultat de régularité
se déduit d'inégalités du type \flqq inégalités de Harnack paraboliques\frqq\ (voir
\cite{Del} et \cite{TLL} pour des applications de cette méthode dans le cas d'espaces discrets).
En utilisant les inégalités gaussiennes du paragraphe précédent,
on est ici plus proche de la démarche de Nash (\cite{N}).

\v1
\n - Le théorème 1.3 se déduit sans difficulté du théorème 3.3,
du théorème 4.1 et des inégalités (\ref{major}) et (\ref{ineg2}), en utilisant l'inégalité triangulaire:
$$
h_{2n} (x)
\leq
\left| h_{2n} (x) - h_{2n} (0) \right|
+
h_{2n} (0)
=
\left| h_{2n} (x) - h_{2n} (0) \right|
+
\norm{h_n}{}^2
\; .
$$ 

\par}

\v5
\noindent{\bf Preuve du théorème 1.2 }
On procède comme dans \cite{TLL}.

\v1
Commençons par remarquer qu'en appliquant le théorème ergodique ponctuel de Birkhoff à $\o{c}$ et à une différence de sommes ergodiques suivant l'idée de la preuve du théorème précédent, on montre que,
pour presque tout environnement, pour tout réel $x_0$, pour tout $\delta>0$, on a:
\begin{eqnarray}\label{conv_erg}
\lim_{n\to +\infty}
{{\o{c} (B (g_n (x_0),\delta\sqrt{2n})\cap 2\Z)}
\over
{\delta\sqrt{2n}}}
=\int \o{c} \; d\mu
\end{eqnarray}

\v1
Par ailleurs, pour tout $\delta>0$, pour tout $n\geq 1$, il vient:
\begin{eqnarray*}
& &
\sqrt{2n}\, {{\P_0 [S_{2n}=g_n (x_0)]} \over{\o{c} (g_n (x_0))}}\, \int \o{c}  \, d\mu
-
2 k_\sigma (x_0)\\
&=&
\sqrt{2n}\, h_{2n} (g_n (x_0)) \, \int \o{c}  \, d\mu
-
{{\o{c} (B (g_n (x_0),\delta\sqrt{2n})\cap 2\Z)}\over{\delta}}\, h_{2n} (g_n (x_0))
\\
& &+\ 
{{1}\over{\delta}}
\sum_{x \in B (g_n (x_0),\delta\sqrt{2n})\cap 2\Z}
h_{2n} (g_n (x_0))\, \o{c} (x)
-
{{1}\over{\delta}}
\sum_{x \in B (g_n (x_0),\delta\sqrt{2n})\cap 2\Z}
h_{2n} (x)\, \o{c} (x)
\\
& &+\ 
{{1}\over{\delta}}
\,
\P_0 [S_{2n} \in B (g_n (x_0),\delta\sqrt{2n})]
-
{{1}\over{\delta}}
\,
\int_{x_0-\delta}^{x_0+\delta} k_\sigma (z)\, dz
\\
& &+\ 
{{1}\over{\delta}}
\,
\int_{x_0-\delta}^{x_0+\delta} k_\sigma (z)\, dz
-
2\, k_\sigma (x_0)
\; ,
\end{eqnarray*}
où l'on a repris les notations du théorème 1.4.

\v2
Ainsi,
\begin{eqnarray}
& &
\left|
\sqrt{2n}\, {{\P_0 [S_{2n}=g_n (x_0)]} \over{\o{c} (g_n (x_0))}}\, \int \o{c}  \, d\mu
-
2 k_\sigma (x_0)\right| \nonumber
\\
&\leq&
\sqrt{2n}\, h_{2n} (g_n (x_0))
\left|
\int \o{c}  \, d\mu
-
{{\o{c} (B (g_n (x_0),\delta\sqrt{2n})\cap 2\Z)}\over{\delta \, \sqrt{2n}}}
\right| \nonumber
\\
& &+\ 
\sqrt{n}
\,
\max_{x\in  B (g_n (x_0),\delta\sqrt{2n})\cap 2\Z}
\left(\left| h_{2n} (g_n (x_0)) - h_{2n} (x) \right|\right)
\,
{{\o{c} \left( B (g_n (x_0), \delta\sqrt{2n})\cap 2\Z) \right)}
\over
{\delta\sqrt{n}}}
\label{equation3}
\\
& &+\ 
{{1}\over{\delta}}
\left|
\P_0 \left[ {{g_n (x_0)}\over{\sqrt{2n}}} - \delta \leq {{S_{2n}}\over{\sqrt{2n}}}
\leq
{{g_n (x_0)}\over{\sqrt{2n}}} + \delta \right]
-
\int_{x_0-\delta}^{x_0+\delta} k_\sigma (z)\, dz
\right| \nonumber
\\
& &+\ 
2
\,
\left(
{{1}\over{2\delta}}
\,
\int_{x_0-\delta}^{x_0+\delta} |k_\sigma (z)-k_\sigma (x_0)|\, dz
\right)\label{equation4}
\; .
\end{eqnarray}

Pour presque tout environnement, un réel $\vareps >0$ étant donné,
on peut fixer $\delta>0$ assez petit de manière à majorer pour tout $n$ assez grand les expressions
$(\ref{equation3})$ et
$(\ref{equation4})$
par $\vareps$. Pour cela, on utilise le théorème 4.1,
la convergence (\ref{conv_erg}) et la continuité de $z\mapsto k_\sigma (z)$ en $x_0$.
On conclut alors en passant à la limite supérieure en $n$
et en utilisant à nouveau la convergence (\ref{conv_erg}),
le théorème 1.3
(pour justifier que la suite $(\sqrt{2n} \, h_{2n} (g_n (x_0)))_{n\geq 0}$ est bornée)
et le théorème limite central associé à la continuité de la loi normale.

\fdem

\v5
\noindent{\bf 5. Une chaîne de Markov analogue en temps continu}

\v2
Dans ce paragraphe, on s'intéresse à un processus stochastique
en temps continu analogue à la marche aléatoire $(S_n)_{n\geq 0}$.
Ce processus est parfois appelé {\it marche aléatoire à vitesse variable}
dans la littérature (voir \cite{TLC6} par exemple).
Parmi les nombreuses manières possibles de l'introduire,
nous choisissons ici celle consistant à partir de $(S_n)_{n\geq 0}$.

\v3
On considère donc à nouveau
un environnement de conductances fixé et
la marche aléatoire $(S_n)_{n\geq 0}$ associée.
Pour chaque réalisation de $(S_n)_{n \geq 0}$, on considère également la réalisation d'une suite
$(T_n)_{n\geq 1}$ de variables aléatoires indépendantes de lois exponentielles
telles que, pour tout $n\geq 1$, la moyenne de $T_n$ soit égale à $\o{c} (S_{n-1})^{-1}$. 

\v2
On introduit alors les {\it instants de sauts} en posant
$J_0:=0$
et, pour tout $n\geq 1$,
$$  
J_n:=T_1+T_2+\cdots+T_n\; .
$$

On pose enfin
$X_t:=S_{n}$,
pour tout $t\geq 0$ satisfaisant $J_n \leq t < J_{n+1}$.

\v1
A environnement de conductances fixé,
on a ainsi défini
un processus stochastique $(X_t)_{t\geq 0}$ qui est en fait une chaîne de Markov à temps continu sur $\Z$.
La marche aléatoire $(S_n)_{n\geq 0}$ étant irréductible et récurrente,
la chaîne de Markov $(X_t)_{t\geq 0}$ est irréductible, non explosive et récurrente.

\v5
Dans la suite, on réunit les éléments permettant de déduire
un théorème limite central local pour $(X_t)_{t\geq 0}$
des méthodes qui ont été développées en temps discret aux paragraphes 3 et 4 ci-dessus.

\v1
Compte tenu des moyennes des temps d'attente $T_n$ en chaque site,
on comprend que la chaîne de Markov $(X_t)_{t\geq 0}$
a \flqq moins de raisons\frqq\ que $(S_n)_{n\geq 0}$ d'être \flqq piégée\frqq\ par les sites qui correspondent
à de grandes valeurs de $\o{c}$. Ceci se traduit en particulier par l'absence
de condition d'intégralité sur $\o{c}$ dans le théorème limite central pour $(X_t)_{t\geq 0}$
que l'on rappelle maintenant.

\v2
\n{\bf Théorème 5.1 }{\it (\cite{K-K}, \cite{D-D}, \cite{L})}

\v1
{\leftskip=1cm
\n (1) Si $1/c$ est intégrable alors,
pour presque tout $\omega$ dans $\Omega$,
pour tous $-\infty \leq a<b\leq +\infty$, on a
$$
\lim_{t\to +\infty}
\P_0^\omega \left[ a<{{X_t}\over{\sqrt{t}}} \leq b\right]
=
\int_{a}^b
k_\sigma (z)
\, dz
\, ,
$$
avec
$$
k_\sigma (z)
:=
{1 \over{\sigma\sqrt{2\pi}}}
\,
\exp \left(-{{z^2}\over{2\sigma^2}}\right)
\quad
{\rm et}
\quad
\sigma^2:={2\over{\int {1\over c}\, d\mu}} \; .
$$

\v1
\n (2) Si $1/c$ n'est pas intégrable alors,
pour presque tout $\omega$ dans $\Omega$,
pour tous $-\infty \leq a<b\leq +\infty$ avec $a$ et $b$ non nuls, on a
$$
\lim_{t\to +\infty}
\P_0^\omega \left[ a<{{X_t}\over{\sqrt{t}}} \leq b\right]
=
\left\{
\begin{array}{cl}
1 & \quad \mbox{si $0\in ]a,b]$}\\
0 & \quad \mbox{sinon}
\end{array}
\right.
\; .
$$
\par}

\v5
D'autre part, la mesure $(\o{c} (x))_{x\in \Z}$ sur $\Z$ étant réversible pour la marche aléatoire $(S_n)_{n\geq 0}$,
la mesure de comptage sur $\Z$ est réversible pour la chaîne de Markov $(X_t)_{t\geq 0}$. Autrement dit,
pour tout $t\geq 0$,
$$
\forall x,y \in \Z\; ,\quad \P_x [X_t=y]=\P_y [X_t=x]\; .
$$
En particulier, $(X_t)_{t\geq 0}$ admet une mesure invariante de masse totale infinie,
elle est donc récurrente nulle et l'on a:
$$
\forall x,y \in \Z\; ,\quad \lim_{t\to +\infty} \P_x [X_t=y]=0 \; .
$$

De plus, le semi-groupe $(P(t))_{t\geq 0}$ associé à $(X_t)_{t\geq 0}$
est une famille d'opérateurs symétriques de $\ell^2 (\Z)$ de norme égale à 1.
Il en résulte que si $\LR$ désigne le générateur associé dans $\ell^2 (\Z)$,
l'opérateur (non borné) $-\LR$ est symétrique
et positif. En outre,
pour toute fonction $f$ à support fini sur $\Z$, on a:
$$
\forall x\in \Z\; ,\quad
\LR f(x)
=
c(x,x+1) f(x+1)+ c(x,x-1)  f(x-1) - \o{c} (x) f(x)\; .
$$

\v3
Le noyau de la chaleur associé à $(X_t)_{t\geq 0}$ est défini par:
$$
h_t (x):=\P_0 [X_t=x]=P(t) \I_{\{0\}} (x) \; ,\quad x\in \Z\; ,\ t\geq 0\; .
$$
Pour tout $t\geq 0$, la symétrie de l'opérateur $P(t)$ donne:
$
(h_t,h_t)=h_{2t} (0)
$
(on a noté $(\cdot , \cdot)$ le produit scalaire dans $\ell^2 (\Z)$).

\v3
Si à présent on pose, pour tout $t\geq 0$, $u(t):=(h_t,h_t)$, on obtient,
en utilisant les équations de Kolmogorov rétrograde et progressive, que
$$
\forall n\geq 0\; ,\quad
u^{(n)} (t)=2^n (h_t,\LR^n h_t)\; ,
$$
où $u^{(n)}$ désigne la dérivée $n$-ième de la fonction $u$.
Ainsi,
$$
\forall n\geq 0\; ,
\quad
(-1)^n u^{(n)} (t)=2^n (h_t,(-\LR)^n h_t) \geq 0
$$
(distinguer les cas $n$ pair et $n$ impair, et utiliser la positivité et la symétrie de $-\LR$).

\v1
Le théorème de Bernstein de représentation des fonctions
complètement monotones sur $[0,+\infty[$ (voir par exemple \cite{F} chapitre XIII)
assure alors l'existence d'une mesure de probabilité $\nu$ sur $[0,+\infty[$
telle que
$$
\forall t\geq 0\; ,
\quad
u(t)=u(0) \int_0^{+\infty} \e^{-tx} \, d\nu (x)
\; .
$$
En particulier, pour tout $t>0$, on a:
\begin{eqnarray*}
2 (h_{2t},-\LR h_{2t})
=
-u'(2t)
&=&
u(0) \int_0^{+\infty} x \e^{-2tx} \, d\nu (x)
\\
&\leq&
\max \{ x\e^{-tx}:\ x\in [0,+\infty[\}\,  u(t)
\\
&=&
{{\e^{-1}}\over{t}}\, u(t)
\; ,
\end{eqnarray*}
ce qui constitue un analogue du corollaire 2.4 ci-dessus.

\v3
Remarquons enfin que, pour toute fonction $f$ appartenant au domaine de $\LR$,
\begin{eqnarray*}
(f,-\LR f)&=&
\lim_{t\to 0^+} \left( f, {1 \over t} \left(f-P(t)f\right)\right)
\\
&=&
\lim_{t\to 0^+} {1\over{2t}} \sum_{x,y\in \Z} (f(x)-f(y))^2 \P_x [X_t=y]
\\
&\geq&
\lim_{t\to 0^+} {1\over 2} \sum_{x,y\in F} (f(x)-f(y))^2 {{\P_x [X_t=y]}\over{t}}
\\
&=&
\sum_{x\in F\, :\ x+1\in F} (f(x)-f(x+1))^2 c(x,x+1)\; ,
\end{eqnarray*}
pour toute partie finie $F$ de $\Z$.

\v5
En utilisant les méthodes des deux paragraphes précédents,
on montre en particulier le théorème suivant.

\v2
\n{\bf Théorème 5.2 }

\v1
{\leftskip=1cm
\n (1) Si $1/c$ est intégrable alors,
pour presque tout environnement de conductances,
pour tout réel $x_0$,
$$
\lim_{t \to +\infty}
\sqrt{t} \, \P_0 \left[ X_t=[x_0 \sqrt{t}] \right]
=
{{\sqrt{\int {1 \over c}\, d\mu}}\over{2\sqrt{\pi}}}
\,
\exp\left( -{{\int {1\over c}\, d\mu}\over{4}}\, x_0^2\right)
$$
(on a noté $[x]$ la partie entière de $x$).

\v2
\n (2) Si $1/c$ n'est pas intégrable alors,
pour presque tout environnement de conductances,
pour tout $x$ dans $\Z$,
$$
\lim_{t \to +\infty}
\sqrt{t} \, \P_0 \left[ X_t=x\right]
=
+\infty
\; .
$$
\par}

\v3
\n{\it Remarques }

\v1
{\leftskip=1cm
\n - La première partie de ce théorème est donnée
sans démonstration dans \cite{K-K}.

\v1
\n - La même méthode s'applique également {\it à la marche aléatoire à vitesse constante},
c'est-à-dire dans le cas où les temps d'attente $T_n$ en chaque site suivent tous
la loi exponentielle de paramètre 1.
\par}

\v9
\n{\bf Remerciements -- }L'auteur remercie Jérôme Depauw, Yves Derriennic
ainsi que le rapporteur 
pour leurs précieuses remarques.

\v5
{\footnotesize

{\footnotesize

\par}

\end{document}